\documentclass[12pt]{amsart}

\usepackage{amsfonts,amsmath,amssymb,amsthm}

\newcommand{\Cc}{\mathbb{C}} 
\newcommand{\Pp}{\mathbb{P}}
\newcommand{\Rr}{\mathbb{R}}
\newcommand{\Nn}{\mathbb{N}}
\newcommand{\Ss}{\mathbb{S}}
\newcommand{\Zz}{\mathbb{Z}}
\newcommand{\defi}[1]{\emph{#1}}

\renewcommand{\epsilon}{\varepsilon}

\renewcommand {\leq}{\leqslant}
\renewcommand {\geq}{\geqslant}

\renewcommand {\Re}{\mathop{\mathrm{Re}}\nolimits}
\renewcommand {\Im}{\mathop{\mathrm{Im}}\nolimits}

\newcommand{\grad}{\mathop{\mathrm{grad}}\nolimits}

{\theoremstyle{plain}
\newtheorem{theorem}{Theorem}[section]    
\newtheorem{lemma}[theorem]{Lemma}       
\newtheorem{proposition}[theorem]{Proposition}      
\newtheorem{corollary}[theorem]{Corollary}      

}
{\theoremstyle{remark}
\newtheorem{definition}[theorem]{Definition}      
\newtheorem*{remark*}{Remark}  
\newtheorem{remark}[theorem]{Remark}   
\newtheorem{example}[theorem]{Example}
}

\begin{document}

\title[Links of meromorphic functions]{Meromorphic functions, bifurcation sets and fibred links}

\author{Arnaud Bodin}
\address{Laboratoire Paul Painlev\'e, Math\'ematiques,
Universit\'e de Lille 1,  59655 Villeneuve d'Ascq, France.}
\email{Arnaud.Bodin@math.univ-lille1.fr}

\author{Anne Pichon}
\address{Institut de Math\'ematiques de Luminy, Campus de Luminy, Case 907
13288 Marseille, France.}
\email{pichon@iml.univ-mrs.fr}
\date{\today}

\begin{abstract}
We give a necessary condition for a meromorphic function in several variables 
to give rise to a Milnor fibration 
of the local link (respectively of the link at infinity).
In the case of two variables we give some necessary and sufficient conditions 
for the local link (respectively the link at infinity) to be fibred.
\end{abstract}

\maketitle

\section{Introduction}

A famous result of J.~Milnor \cite{M} states that the link
$f^{-1}(0) \cap {\Ss}^{2n-1}_\epsilon$ ($0<\epsilon \ll 1$) of a 
  holomorphic germ $f : (\Cc^n,0) \longrightarrow (\Cc,0)$  is a fibred
link and moreover that a fibration is given by the so-called Milnor
fibration $\frac{f}{|f|} : {\mathbb S}^{2n-1}_\epsilon \setminus f^{-1}(0)
\longrightarrow \Ss^1.$
Throughout this paper ${\Ss}_{r}^{n-1}$ denotes the sphere with radius $r$ centered at the origin of 
$\Rr^{n}$. 

The proof of this result has been extended in several directions in order to construct some natural fibrations in other situations of singularity theory. In this paper, we focus on two of them : 
\begin{enumerate}
\item  Let $f : ({\mathbb
 R}^{n+k},0) \longrightarrow ({\mathbb
 R}^k,0)$ be a real analytic  germ  with an isolated critical point at the origin. J.~Milnor  \cite[Chapter 11]{M} proved that  for every sufficiently small sphere ${\mathbb
 S}^{n+k-1}_{\epsilon} $ centered at the origin  in $ {\mathbb
 R}^{n+k}$, 
the complement $ {\mathbb S}^{n+k-1}_{\epsilon} \setminus L_f$ of the link $L_f =  {\mathbb
 S}^{n+k-1}_{\epsilon}  \cap f^{-1}(0)$  fibres over the circle.  As pointed out by Milnor,  the fibration is not
necessarily given by the Milnor map $\frac{f}{\|f\|}$.  

This result can be extended  to a real analytic germ $f:(X,p) \longrightarrow ({\mathbb
 R}^k,0) $   with isolated critical value satisfying a suitable stratification condition, where $(X,p)$ is  a germ of real analytic space with isolated singularity at $p$ (\cite{PS}, Theorem 1.1).

\item Another direction deals with links at infinity. Let $f : \Cc^n \longrightarrow \Cc$ be a polynomial map. The link  at infinity associated with the fibre $f^{-1}(0)$ is defined by 
$L_{f,\infty}=f^{-1}(0) \cap {\mathbb S}^{2n-1}_R$ for a sufficiently large radius $R \gg
1$. In \cite{NZ}, A.~N\'emethi and A.~Zaharia
 proved that under a
condition called \emph{semitame},  the link at infinity $L_{f,\infty}$ is fibred by 
 the Milnor map $\frac{f}{|f|} : \Ss^{2n-1}_R \setminus L_{f,\infty} \longrightarrow {\mathbb S}^1$. 
 
 When $n=2$, A.~Bodin in \cite{Bo} proved that
the link at infinity is fibred if and only if the semitame condition
holds, or equivalently iff the set of critical values at infinity is
void or equal to $\{0\}$.
\end{enumerate}

In the first part of the paper we consider a meromorphic function $f/g$, where $f : (\Cc^n,0) \longrightarrow (\Cc,0)$ and $g: (\Cc^n,0) \longrightarrow (\Cc,0)$ are two holomorphic germs. We firstly define a condition 
called \defi{semitame at the origin} for $f/g$ by adapting the definition of \cite{NZ}.  Namely it consists of defining a bifurcation set $B$ which reflects the behaviour of the  points of non-transversality between the fibres of $f/g$ and the spheres ${\Ss}^{2n-1}_{\epsilon}$, $\epsilon \ll 1$,  centered at the origin of $\mathbb C^n$. In particular, under the semitame condition $B=\{0,\infty\}$, there is no non-isolated singular points of the meromorphic function $f/g$ (nor point of indeterminacy) outside $(f=0) \cup (g=0)$.

\begin{theorem}
\label{th:intro1}
Let  $f,g :  (\Cc^n,0) \longrightarrow (\Cc,0)$ be two germs of holomorphic functions without common factors such that the meromorphic germ $\frac fg$  verifies the semitame condition $B=\{0,\infty\}$
at the
origin.  Then for all sufficiently small $\epsilon >0$ the Milnor map 
$$\frac{f/g}{|f/g|} : \Ss^{2n-1}_\epsilon \setminus (L_f \cup L_g) \longrightarrow \Ss^1$$
is a $C^{\infty}$ locally trivial fibration.  

Moreover for $n=2$ this fibration is a fibration of the multilink $L_{f/g} = L_f \cup - L_g$.
\end{theorem}

\bigskip

The second part is devoted to find a reciprocal of Theorem \ref{th:intro1}
in the case of two variables: $n=2$. More precisely we give other equivalent conditions
to the equivalence (\ref{it:i1}) $\Leftrightarrow$ (\ref{it:i2}) due to A.~Pichon and J.~Seade.
\begin{theorem}
\label{th:intro2}
Let $f: ({\mathbb C}^2,0) \to ({\mathbb C},0)$ and $g: ({\mathbb C}^2,0) \to ({\mathbb C},0)$ be two holomorphic germs without common branches. Then the following conditions are equivalent:
\begin{enumerate}
\item \label{it:i1} The multilink $L_f  \cup -L_g$ is fibred;
\item  \label{it:i2} The real analytic germ $f \overline g : ({\mathbb C}^2,0) \to ({\mathbb C},0)$ has $0$ as an isolated critical value;
\item \label{it:i3} The Milnor map  $\frac{f \overline{g} }{ |{f \overline {g}}|} : 
 {\mathbb S}^{3}_{\epsilon} \setminus  (L_f \cup   L_g)$  is a $C^{\infty}$ locally trivial fibration;
 \item \label{it:i4} The meromorphic map $f/g$ holds the semitame condition at the origin;
 \item \label{it:i5} Each $c \neq 0, \infty$ is a generic value of the local pencil generated by $f$ and $g$.
\end{enumerate}
\end{theorem}

The equivalence (\ref{it:i1}) $\Leftrightarrow$ (\ref{it:i2}) is proved by A. Pichon in \cite{P} when $f$ and $g$ have $0$ as an isolated critical point and generalized by A. Pichon and J. Seade in \cite{PS}.
In particular Theorem \ref{th:intro2} gives an alternative proof to $(1) \Rightarrow (5)$ which was first observed by F.~Michel and H.~Maugendre in \cite{MM}.

Notice that a natural question consists  of  comparing the Milnor fibration $\frac{f \overline{g} }{ |{f \overline {g}}|} :  {\Ss}^{3}_{\epsilon} \setminus  (L_f \cup   L_g) \longrightarrow \Ss^1$ with the local fibrations of meromorphic germs introduced in \cite {GLM1, GLM2}; this will be done in a forthcoming work with J.~Seade. 

\bigskip

The last part of the paper is devoted to the global case. Let 
$f: {\mathbb C}^2 \to {\mathbb C}$ and $g: {\mathbb C}^2 \to {\mathbb C}$ be two polynomial maps. The link at infinity $L_{f/g,\infty}$ of the meromorphic function $f/g$ is defined by $L_{f/g,\infty} = L_{f,\infty} \cup -L_{g,\infty}$. In \cite{HR}, after observing that $\frac{f \bar g}{| f \bar g |} = \frac{f/g}{|f/g|}$, M.~ Hirasawa and L.~Rudolph ask whether the methods developed in \cite{P} can be adapted to  $f/g$ at infinity. Namely there are two natural questions:  Under which conditions is the link $ L_{f/g,\infty}$  fibred?  When $ L_{f/g,\infty}$  is fibred, is the Milnor map $\frac{f/g}{|f/g|} : {\mathbb S^3_{R}} \setminus L_{f/g,\infty} \to {\mathbb S}^1$ a fibration of  $L_{f/g,\infty}$?

We define a semitame condition at infinity for meromorphic maps, that enables one to
adapt the methods of \cite{NZ} and \cite{Bo} and of the first sections of this work. We first 
get a version of Theorem \ref{th:intro1}  at infinity (see Theorem \ref{th:fibinfty}) and then we obtain the following result, which is a complete answer to the question of M.~Hirasawa and L.~Rudolph: 

\begin{theorem}
\label{th:intro3}
  Let $f: {\mathbb C}^2 \to {\mathbb C}$ and $g: {\mathbb C}^2 \to {\mathbb C}$ be two polynomial maps. The following conditions are equivalent: 
\begin{enumerate}
\item the meromorphic function $f/g$ is semitame at infinity;
\item the link $L_{f/g,\infty}$ is fibred.
\end{enumerate}
Moreover, if these conditions hold, then the Milnor map $\frac{f/g}{|f/g|} : {\mathbb S^3_{R}} \setminus L_{f/g,\infty} \to {\mathbb S}^1$ is a fibration of the link $L_{f/g,\infty}$.
\end{theorem}

It should be noticed that the local situation for meromorphic maps is
similar to the polynomial situation at infinity and the situation at
infinity for meromorphic maps is morally the gluing of a finite number
of meromorphic local situations.

\bigskip

\emph{Acknowledgements:} It is a pleasure to thank Lee Rudolph whose thoughts and
questions
are a perpetual source of inspiration.

\section{Milnor map of a meromorphic function:\\ the local case}
\label{sec:loc}

\subsection*{The multilink of a meromorphic function} 

An {\it oriented link} in $\mathbb S^k$ is a disjoint union of
oriented $(k-2)$-spheres $K_1,\ldots,K_\ell$ embedded in $\mathbb S^k$. A
multilink is the data of an oriented link $L=K_1 \cup \ldots \cup K_\ell$
together with a multiplicity $n_i \in {\mathbb Z}$ associated to
each component $K_i$ of $L$. We denote $L=n_1 K_1 \cup \ldots \cup n_\ell
K_\ell$ with the convention $n_i K_i = (-n_i) (-K_i)$ where $-K_i$ means
$K_i$ with the opposite orientation (see \cite{EN}).

For example, the multilink of a holomorphic germ $f : ({\mathbb
  C}^n,0) \to ({\mathbb C},0)$ is defined by $L_f = n_1 L_{f_1} \cup
\ldots \cup n_\ell L_{f_\ell}$, where $f = \prod _{i=1}^\ell f_i^{n_i}$ is
the decomposition of $f$ in the \textsc{ufd} ${\mathbb C}\{x_1,\ldots , x_n
\}$ and where $L_{f_i} = f_i^{-1}(0) \cap {\mathbb
  S}^{2n-1}_{\epsilon}$, $\epsilon \ll 1$, is oriented as the boundary of
the piece of complex manifold $f_i^{-1}(0) \cap {\mathbb
  B}^{2n}_{\epsilon}$. The diffeomorphism class of the pair
$(\Ss^{2n-1}_\epsilon,L_f)$ is independent of $\epsilon$ when
$\epsilon$ is sufficiently small.

\begin{definition} 
  Let  $f,g :  (\Cc^n,0) \longrightarrow (\Cc,0)$ be two germs of holomorphic
  functions without common factors in the decomposition in the \textsc{ufd}
  ring $\Cc\{x_1,\ldots, x_n\}$. Let $f,g : U \to \Cc$ be some
  representative of $f$ and $g$.  The \defi{multilink} $L_{f/g}$ of
  the meromorphic function $f/g : U \longrightarrow \Pp^1$ is defined by
    $$L_{f/g} = L_f \cup -L_g$$
  where $L_f$ and $L_g$ denote the
  multilink of $f$ and $g$.
\end{definition}

\subsection*{Semitame map at the origin}

Let  $U$ be an open neighbourhood of $0$ in $\Cc^n$ and let $f,g :  U \longrightarrow \Cc$ be two holomorphic functions without common factors such that $f(0) = g(0)=0$. 

Let us consider the meromorphic function $f/g : U \rightarrow {\Pp}^1 $ defined by $(f/g) (x) =  [f(x) : g(x)]$. Notice that $f/g$ is not defined on the whole $U$; its indetermination locus is $I(f/g) = \{ x \in U \ | \ f(x)=0 \text{ and } g(x)=0 \} $. 

Adapting Milnor's definition let us define the
gradient of $f/g$ outside $I(f/g)$ by:
$$\grad (f/g) = \left(\overline{\frac{\partial (f/g)}{\partial x_1}},
  \ldots, \overline{\frac{\partial (f/g)}{\partial x_n}}\right).$$

Let us consider the set
$$M(f/g) = \left\lbrace x \in U \setminus I(f/g) \mid \exists \lambda \in \Cc,
  \grad \frac fg (x)=\lambda x \right\rbrace$$
consisting of the points of
non-transversality between the fibres of $f/g$ and the spheres ${\mathbb
  S}^{2n-1}_{r}$ centered at the origin of $\Cc^n$.

We define a bifurcation set $B \subset \Pp^1$ for the
meromorphic function $f/g$ as follows.
\begin{definition}\label{semi0}
  The {bifurcation set} $B$ consists of all values $c\in \Pp^1$ such that there
  exists a sequence $(x_k)_{k \in {\mathbb N}}$ of points of $M(f/g)$
  such that
  $$\lim_{k \to \infty} x_k = 0 \quad \text{and} \quad \lim_{k \to \infty}
  \frac{f(x_k)} {g(x_k)}= c.$$
  By convention and to avoid discussion
  we set
  $$\{0,\infty\}\subset B.$$
\end{definition}

The following definition is adapted from that of \cite{NZ} which concerned
polynomial maps at infinity.
\begin{definition}\label{semitame}
  The meromorphic germ $f/g$ is \defi{semitame at the origin} if $B =
  \{0,\infty\}$.
\end{definition}

\begin{remark} \ 
\begin{enumerate}
\item  When $f/g$ is semitame at the origin, the non isolated singular points of the meromorphic function $f/g$ in an open  neighbourhood of the origin belong to $(f=0) \cup (g=0)$. Indeed, assume that there exists a sequence $(x_k)_{k \in {\mathbb N}}$ of points of $U \setminus I(f/g)$ such that 

$$\lim_{k \rightarrow \infty} x_k = 0  \hbox{ and } \grad (f/g)(x_k) = 0.$$ 

Then for all $k$, $x_k \in M(f/g) $. As the critical values of $f/g$ are isolated, one can assume that there exists $c \in {\Pp}^1 $ such that for all $k$, $(f/g)(x_k) = c$. The semitame condition therefore  implies $c=0$ or $c=\infty$. 

\item Notice that $f$ and $g$ can have non-isolated singularities,
  whereas $\frac{f}{g}$ is semitame at the origin. See example
  \ref{ex}.
  
\item It is not hard to prove that the bifurcation set of $g/f$ is the
  set of the inverse elements $\frac{1}{c}$ of the elements $c$ of the
  bifurcation set of $f/g$.
  
\item One can prove (using e.g. the arguments of the proof of Lemma
  \ref{lem:arg}) that a sequence $(x_k)$ as in Definition
  \ref{semitame} verifies:
  $$\lim_{k \to \infty} \|x_k\| \cdot \|\grad \frac fg (x_k)\| = 0.$$
\end{enumerate}
\end{remark}

\begin{example}
\label{ex}
\begin{enumerate}
\item  Let $f(x,y) = x^2$ and $g(x,y) = {y^3}$. Then $B=\{0,\infty\}$
  and $f/g$ is semitame at the origin.
  
\item Let $f(x,y) = x(x+y^2)+y^3$ and $g(x,y) = {y^3}$. Then
  $B=\{0,1,\infty\}$ and $f/g$ is not semitame at the origin.
\end{enumerate}
\end{example}

\subsection*{Fibration theorem under the semitame condition}

\begin{theorem}
\label{th:fib}  
Let  $f,g :  (\Cc^n,0) \longrightarrow (\Cc,0)$ be two germs of
 holomorphic functions without common factors such that the meromorphic 
germ $\frac fg$  is semitame at the
origin.  Then there exists $0<\epsilon_0 \ll 1$ such that for each
$\epsilon \leq \epsilon_0$ the Milnor map
$$\frac{f/g}{|f/g|} : \Ss^{2n-1}_\epsilon \setminus L_{f/g}
\longrightarrow {\mathbb S}^1$$
is a $C^{\infty}$ locally trivial fibration.
\end{theorem}

The proof follows Milnor's proof \cite[Chapter 4]{M} with minor
modifications.  See also \cite{NZ}.
The main modification concerns Lemma 4.4 of \cite{M}, for which we
give an adapted formulation and a detailed proof:

\begin{lemma}
\label{lem:arg} Assume that the meromorphic germ $G=f/g$ is semitame 
at the origin. Let $p:[0,1] \longrightarrow \Cc^n$
be a real analytic path with $p(0)=0$ such that for all $t>0$, $G(p(t))
\notin \{0,\infty\}$ and such that the vector $\grad \log G(p(t))$ is
a complex multiple $\lambda(t)p(t)$ of $p(t)$. Then the argument of the complex
number $\lambda(t)$ tends to $0$ or $\pi$ as $t\rightarrow 0$.
\end{lemma}

\begin{proof}
  The equality $\grad \log G(p(t)) = \lambda(t)p(t)$ implies $\grad
  G(p(t)) = \lambda(t) p(t) \overline{G(p(t))}$. Therefore $p(t) \in
  M(G)$.
  
  Let us consider the expansions
\begin{eqnarray*}
 p(t) = \mathbf {a} t^\alpha+\ldots,   \\
 G(p(t)) = bt^\beta+\ldots, \\
 \grad G(p(t)) = \mathbf{c}t^{\gamma}+\ldots,
\end{eqnarray*}
with $\alpha \in \Nn^*$, $\beta, \gamma \in \Zz$, $\mathbf {a} \neq
0$, $b\neq 0$.

Assume that $\mathbf{c}=0$, {\it i.e.} that $\grad G(p(t))$ is
identically $0$. By definition of $\grad(G)$, one has \begin{equation}
\label{eq:dG}
\forall t \in ]0,1[, \ \ \ \frac{dG}{dt}(p(t))= \left\langle 
\frac{dp}{dt}\mid \grad G (p(t)) \right\rangle.
\end{equation}
Therefore $G(p(t))$ is a constant $\nu$, different from $0$ and
$\infty$ by the hypothesis of the lemma. This contradicts the fact
that $G$ is semitame as $\nu$ belongs to the bifurcation set $B$. Then
in fact $\mathbf{c} \neq 0$.

Replacing each term by its expansion in the equality $\grad G(p(t)) =
\lambda(t) p(t) \overline{G(p(t))}$, we get
$$
{\mathbf {c}}t^{\gamma}+\ldots=\lambda(t)(\mathbf {a}
t^\alpha+\ldots)(\bar{b}t^{\beta}+\ldots).$$

Identifying the coefficient of lower degree, we obtain $\lambda(t) =
\lambda_0t^{\gamma-\alpha-\beta}+\ldots$, and $\mathbf {c}=
\lambda_0 \mathbf {a}\bar b$.  Then from (\ref{eq:dG}) we obtain
$$\beta b t^{\beta-1}+\ldots = \alpha \|a\|^2\bar\lambda_
0bt^{\alpha-1+\gamma}+\ldots$$

Assume that $\beta=0$ then $\lim_{t\rightarrow 0} G(p(t)) =b \in
\Cc^*$ (and $\alpha+\gamma>0$ which implies $\|p(t)\|\cdot\|\grad
G(p(t))\| \rightarrow 0$) and $b$ belongs to the bifurcation set $B$.
This contradicts the fact that $G$ is semitame. Then in fact $\beta
\neq 0$.

Therefore, $\beta = \alpha \|a\|^2\bar\lambda_ 0$ which proves that
$\lambda_0$ is a non-zero real number.
\end{proof}

Another modification of the proof of Milnor to apply the
curve selection lemma is to transform all equalities
involving meromorphic function into real analytic equalities.  For
instance, for $G = \frac fg$  let us consider the set $M(G)$ of all $z\in \Cc^n$ for
which
the vectors $\grad G(z)$ and $z$ are linearly dependent, as in Lemma
4.3 of \cite{M}.  Then the equation $\grad G(z) = \lambda z$, where
$z=(z_1,\ldots,z_n)$, is equivalent to the system of analytic
equations:
$$f \frac{\partial g}{\partial z_i}-g \frac{\partial f}{\partial z_i}=
\lambda z_i \bar g^2, \qquad i=1,\ldots,n.$$
As in Milnor's proof
these equations can be transformed into real analytic equations with
real variables $(x_1,y_1,\ldots,x_n,y_n)$, where $x_i =\Re(z_i)$ and $
y_i=\Im(z_i)$. Then the set $M(G)$ is a real
analytic set.

\bigskip

\textbf{In the next paragraphs, we restrict ourselves to the case $n=2$.}

\subsection*{The Milnor map and the local multilink}

At first, let us recall the notion of fibration of a multilink in
${\mathbb S^3}$. For more details see \cite{EN}.

\begin{definition}
  A multilink $L=n_1 K_1 \cup \ldots \cup n_\ell K_\ell$ in ${\mathbb S}^3$
  is fibred if there exists a map $\Phi : {\mathbb S}^3 \setminus L
  \longrightarrow {\mathbb S}^1$ which satisfies the following two
  conditions:
\begin{enumerate}
\item The map $\Phi$ is a $C^{\infty}$ locally trivial fibration;
\item For each $i =1,\ldots,\ell$, let $m_i = \partial D_i$ be the
  boundary of a meridian disk of a small tubular neighbourhood of $K_i$
  oriented in such a way that $D_i . K_i = +1$ in ${\mathbb S}^3$. The
  degree of the restriction of $\Phi$ to $ m_i$ equals $n_i$.
\end{enumerate}
\end{definition}

The following is obtained by examining the behaviour of the map $\pi
\circ \frac{f/g}{|f/g|}$ near each component of the strict transform
of $fg$, where $\pi : \Sigma \to {\mathbb C}^2$ denotes a resolution
of the germ $fg$. For details see \cite[Proposition 3.1]{P} or
 \cite[Lemma 5.1]{PS}.
\begin{proposition}
\label{binding}
  Let $f,g : (\Cc^2,0) \longrightarrow (\Cc,0)$ be two germs of
  holomorphic functions without common factors. If the Milnor map
  $$\frac{f/g}{|f/g|} : \Ss^{3}_\epsilon \setminus L_{f/g}
  \longrightarrow {\mathbb S}^1$$
  is a $C^{\infty}$ locally trivial fibration,
  then it is a fibration of the multilink $L_{f/g}$.
\end{proposition}

\subsection*{The bifurcation set $B$ and the special fibres of the pencil}

Let $f,g : (\Cc^2,0) \longrightarrow (\Cc,0)$ be two germs of
holomorphic functions without common factors.  Let $V \subset \Cc $ be a
Zariski open. Let us denote by $(f/g)^{-1}(t)$ the germ of plane curve
at the origin of $\Cc ^2$ with equation $f(x,y) - tg(x,y)=0$. The
pencil of curves $((f/g)^{-1}(t))_{t \in V}$ is \defi{equisingular} if
for all $t_1,t_2\in V$, the curves $(f/g)^{-1}(t_1)$ and
$(f/g)^{-1}(t_2)$ are equisingular in the sense of Zariski. There
exists a maximal $V_{\text{max}}$ with this property (see e.g. \cite{LW}).
  
\begin{definition}
   The set $B' = \Pp^1 \setminus V_{\text{max}}$ is called the set of
   the \defi{special fibres} of the pencil of plane curves generated
   by $f$ and $g$, or equivalently, the special fibres of the
   meromorphic function $f/g$.
 
   By convention and to avoid discussions we set
   $$\{0,\infty\}\subset B'.$$
\end{definition}

The following result is a meromorphic and local version of parts of the
well-known equivalence of the different definitions of a critical
value at infinity for a polynomial map, see \cite{D} for a survey and proofs, and
also \cite{Ha},\cite{LW}, \cite{Pa}, \cite{Ti}. 
\begin{proposition}
\label{equiv}
$(n=2)$ \\
  For $c\notin\{0,\infty\}$, the following assertions are equivalent:
\begin{enumerate}
\item $c \notin B$;
\item The topological type of the germ of curve $(f/g)^{-1}(t)$ is
  constant for all $t$ near $c$;
\item $c$ is a regular value of the map $\Phi$ obtain from the
  resolution of $f/g$ at the origin;
\item $c \notin B'$.
\end{enumerate}
\end{proposition}

\begin{corollary}$(n=2)$
\label{prop:equality}
$$B = B'.$$
\end{corollary}

\subsection*{A condition for fibration in terms of bifurcation sets} 

\begin{proposition}
\label{th:nonfib}
Let $f,g : (\Cc^2,0) \longrightarrow (\Cc,0)$ be two germs of
holomorphic functions without common factors.  If the link $L_{f/g}$
of $f/g$ is fibred then $B'=\{0,\infty\}$.
\end{proposition}

\begin{proof}
  The proof follows the one of \cite[Theorem 2]{Bo} excepted for the
  case where the link $L_{f/g}$ (or the underlying link if it has
  multiplicities) is the Hopf link.
  
  We briefly recall the ideas from the proof. We suppose that $c\in
  B'$ with $c\notin\{0,\infty\}$, then
  $$F = \left( \frac {f/g}{|f/g|}\right)^{-1}\left(-\frac{c}{|c|}\right)\cap
  \Ss^3_\epsilon,$$
  is a Seifert surface for the link $L_{f/g}$. According to Proposition \ref{equiv} 
  there exists a dicritical divisor $D$ of the
  resolution of $f/g$ at $0$ that is of valency $3$, that is to say it
  corresponds to a Seifert manifold in the minimal Waldhausen
  decomposition of $\Ss^3_\epsilon \setminus L_{f/g}$. For $\omega$
  sufficiently near $c$, there exists a connected component $\ell$ of
  the link $(f/g)^{-1}(\omega)\cap {\mathbb S}^3_\epsilon$ corresponding to $D$.
  Clearly $\ell \cap F = \varnothing$. Then by the characterisation of
  fibred links in \cite[Theorem 11.2]{EN} if $L_{f/g}$ is not the Hopf
  link then $L_{f/g}$ is not fibred.
  
  For the Hopf link, up to an analytical change of coordinates we can
  suppose $f(x,y)=x^p$ and $g(x,y)=y^q$, that implies
  $B'=\{0,\infty\}$.
\end{proof}

\begin{remark} 
  There is an alternative proof using the results of \cite{Ma}  
  and \cite{MM} about another bifurcation set $B''$ defined in
  terms of the Jacobian curve of the morphism $(f,g) : ({\mathbb
    C}^2,0) \to ({\mathbb C}^2,0)$:  if the multilink
  $L_{f/g}$ is fibred, then \cite[Theorem 1.1]{Ma} implies that $1$
  is not a Jacobian quotient of the germ $(f,g) : ({\mathbb C}^2,0)
  \to ({\mathbb C}^2,0)$. Then, by Remark 3 and Theorem 1 of
  \cite{MM}, any $c \not\in \{0,\infty\}$ is a generic value of the
  pencil generated by $f$ and $g$.
\end{remark}

\subsection*{Summary} 

Theorem \ref{th:fib}, Corollary \ref{prop:equality} and Proposition \ref{th:nonfib} lead to  the following Theorem : 
\begin{theorem}
\label{th:eq}
 Let  $f,g :  (\Cc^2,0) \longrightarrow (\Cc,0)$ be two germs of holomorphic functions without common factors. The following are equivalent:
\begin{enumerate}
\item \label{Gsem} The meromorphic function $f/g$ is semitame at the origin;
\item \label{Gfib} Each $c \not\in \{0,\infty\}$ is a generic value of the pencil of curves generated by $f$ and $g$;
\item \label{Kfib} The multilink $L_{f/g}$ is fibred.
\end{enumerate}

Moreover, if these conditions hold, then the Milnor map 
$$\frac{f/g}{|f/g|} : {\mathbb S}^3_{\epsilon} \setminus L_{f/g} \longrightarrow {\mathbb S}^1$$
is a  fibration of the multilink $L_{f/g}$.
\end{theorem}

\begin{proof}
  (\ref{Gsem}) $\Rightarrow$ (\ref{Kfib}) is Theorem \ref{th:fib} and Proposition \ref{binding} \\
  (\ref{Kfib}) $\Rightarrow$ (\ref{Gfib}) is Proposition \ref {th:nonfib} \\
  (\ref{Kfib}) $\Leftrightarrow$ (\ref{Gsem}) is Proposition \ref{prop:equality}
\end{proof}

Now, Theorem \ref{th:eq} and Theorem 2 of \cite{PS} can be summarized in the  statement of Theorem \ref{th:intro2} of the introduction.


\section{Milnor map of a meromorphic function:\\ the global case}

We now produce a very similar description for singularities at
infinity, which leads to a complete answer to the question of Hirasawa and Rudolph. The statements and the proofs are directly adapted from that of Section 2.

Let $f,g  \in \Cc[x_1,\ldots,x_n]$ be two
polynomials with no common factor.  For short, we denote by $f/g$ the meromorphic map $f/g : {\mathbb C}^2  \longrightarrow {\mathbb C}$ well-defined outside $I(f/g)$. Recall the $L_{f,\infty}$ and $L_{g,\infty}$ denotes the multilinks at infinity $f^{-1}(0) \cap \Ss^{2n-1}_R$ and $g^{-1}(0) \cap \Ss^{2n-1}_R$, $R \gg 1$,  of $f$ and $g$ respectively.
\begin{definition}
The \defi{multilink at infinity} $L_{f/g,\infty}$ of the meromorphic function $f/g$ is defined by 
$$ L_{f/g,\infty} = L_{f,\infty} \cup - L_{g,\infty}$$
\end{definition}

We will state a fibration theorem for $L_{f/g,\infty}$  under
a semitame condition. Then for $n=2$ we will state  the reciprocal.
We again consider  the set $M(f/g)$ defined as in Section \ref{sec:loc} (where $U$ is now $\Cc^n$) and we define a bifurcation set $B_\infty \subset \Pp^1$ for the meromorphic function $f/g$ at infinity as follows. 
\begin{definition}
The set $B_\infty $ consists of all values $c\in \Pp^1$ such that there
  exists a sequence $(x_k)_{k\in \Nn}$ of $M(G)$ such that
  $$\lim_{k \to \infty} \|x_k\| =+\infty \quad \text{and} \quad \lim_{k \to \infty} G(x_k)
  = c.$$
  By convention,  we set
  $$\{0,\infty\}\subset B_\infty.$$
\end{definition}

\begin{definition}
  The meromorphic map $f/g$ is \defi{semitame at infinity}
  if $$B_\infty = \{0,\infty\}.$$
\end{definition}

 As for polynomial maps at infinity (\cite{NZ}, \cite{Bo}) we can state a Milnor fibration theorem which is the adaptation of Theorem \ref{th:fib} and Proposition \ref{binding} to the global case. 
\begin{theorem}
\label{th:fibinfty}
Let $f,g  \in \Cc[x_1,\ldots,x_n]$ be two
polynomials with no common factor. If the meromorphic function  is semitame at infinity, then there exists 
 $R_0 \gg 1$ such that for each $R \geq R_0$, the Milnor map 
$$\frac{f/g}{|f/g|} : {\mathbb S}^{2n-1}_R \setminus L_{f/g,\infty} \longrightarrow {\mathbb S}^1$$
is a $C^{\infty}$ locally trivial fibration.
Moreover it is a fibration of the multilink at infinity  $L_{f/g,\infty}$.

\end{theorem}

In the sequel, we restrict to the case $n=2$. Let $f,g  \in \Cc[x,y]$ be two
polynomials with no common factor. We denote by $\tilde f \in \Cc[x,y,t]$ and $\tilde g \in \Cc[x,y,t]$ the homogenisations of $f$ and $g$. The meromorphic map $\tilde{G} = [\frac{\tilde f}{\tilde g} : 1] :  {\mathbb P}^2 \longrightarrow {\mathbb P}^1$ may not be defined at some point on the line
at infinity $H_\infty = \{t=0\}$ of $\Pp^2$.  In the case $\deg f = \deg g$,
then the restriction $\tilde G_| : H_\infty \longrightarrow \Pp^1$ is
a ramified covering.

\begin{definition}
  A point of $H_\infty$ where $\tilde G$ is not well-defined is an
  \defi{indetermination point} of $\tilde G$, and a point of $H_\infty$  where $\tilde G$ is
  well-defined but the restriction $\tilde G_| : H_\infty \rightarrow
  \Pp^1$ is ramified is a \defi{ramification point} of $\tilde G$.
\end{definition}

Let $\pi : \Sigma \longrightarrow {\mathbb P}^2$ be a resolution of the meromorphic function $\tilde{G}$, {\it i.e.} the composition of a finite sequence of blows-up of points starting with the blows-up of the indetermination and of the ramification points of $\tilde{G}$ such that there exists a map $\hat{G} : \Sigma \longrightarrow {\mathbb P}^2$ such that $\hat{G} = \tilde{G} \circ \pi$ which is well defined on $\pi^{-1}( H_\infty)$

The following is the analogous of Proposition \ref{equiv}. 
\begin{proposition}
\label{prop:equalityinfty}
For $n=2$, $c\notin\{0,\infty\}$, the following assertions are equivalent:
\begin{enumerate}
\item $c \notin B_\infty$;
\item outside a large compact set of $\Cc^2$, the topological type of the curve 
  ${f/g}^{-1}(s)$ is constant for all $s$ near $c$;
\item $c$ is a regular value of the map $\hat{G} $.

\end{enumerate}
\end{proposition}

Proposition \ref{prop:equalityinfty} and the arguments of the proof of Proposition \ref{th:nonfib} (or \cite{Bo}, Theorem 2) lead to the following:

\begin{proposition}
\label{th:nonfibinfty} Let $f,g  \in \Cc[x,y]$ be two
polynomials with no common factor. 
If the multilink at infinity $L_{f/g,\infty}$ of $f/g$ is fibred 
then $B_\infty=\{0,\infty\}$.
\end{proposition}

These results enable one to answer positively to the  question of
M.~Hirasawa and L.~Rudolph \cite{HR}:
\begin{theorem}
\label{th:eqinfty}
Let $f,g  \in \Cc[x,y]$ be two
polynomials with no common factor.  The following conditions are equivalent:
\begin{enumerate}
\item \label{Gseminfty} The meromorphic map $f/g$ is semitame at infinity;
\item \label{Kfibinfty} The multilink at infinity $L_{f/g,\infty}$ is  fibred.
\end{enumerate}

Moreover, if these conditions hold, then the Milnor map 
$$\frac{f/g}{|f/g|} : {\mathbb S}^{3}_R \setminus L_{f/g,\infty} \longrightarrow {\mathbb S}^1$$
is a fibration of  the multilink  $L_{f/g,\infty}$.
\end{theorem}

\begin{example}
  Let us take an example from \cite{HR} :  $G= \frac{x^2-y^2}{x^2+y^2}$. 
  Then the link at infinity is not fibred and we have
  $B_\infty=\{0,+1,-1,\infty\}$.
\end{example}



\begin{thebibliography}{22}
  
\bibitem{Bo} A. Bodin, Milnor fibration and fibred links at
  infinity. Inter. Math. Res. Not. 11 (1999), 615--621.
  
\bibitem{D} A. Durfee, Five definitions of critical point at
  infinity. Singularities (Oberwolfach, 1996). Progr. Math. 162,
  Birkh\u{a}user, Basel (1998), 345--360.
  
\bibitem{EN} D. Eisenbud, W. Neumann, Three-dimensional link
  theory and invariants of plane curve singularities. Annals of
  Mathematics Studies 110, Princeton University Press (1985).
  
\bibitem{GLM1}  S. M. Gusein-Zade, I.  Luengo, A.  Melle-Hern\'andez, 
  Zeta functions of germs of meromorphic functions, 
and the Newton diagram. Funct. Anal. Appl.  32  (1998), 93--99.

\bibitem{GLM2}  S. M. Gusein-Zade, I.  Luengo, A.  Melle-Hern\'andez, 
 On the topology of germs of meromorphic functions and its 
applications,  St. Petersburg Math. J.  11  (2000), 775--780.

\bibitem{Ha} Ha, H. V., Nombres de Lojasiewicz et singularit\'es \`a
  l'infini des polyn\^{o}mes de deux variables complexes. C. R.
    Acad. Sci. Paris S\'er. I Math. 311 (1990), 429--432.
  
\bibitem{HR} M. Hirasawa, L. Rudolph, Constructions of Morse maps
  for knots and links and upper bounds on the Morse-Novikov number.
  Preprint (2003).
  
\bibitem{LW} D. T. L\^e , C. Weber, \'Equisingularit\'e dans les
  pinceaux de germes de courbes planes et $C\sp 0$-suffisance.
  Enseign. Math.  43 (1997), 355--380.
  
\bibitem{Ma} H. Maugendre,  Discriminant d'un germe $(g,f): {\mathbf C}^2 \to {\mathbf C}^2$ 
et quotients de contact dans la r\'esolution de $f.g$, Annales de la Facult\'e des 
Sciences de Toulouse, 7 (1998), 497--525.

\bibitem{MM} H. Maugendre, F. Michel, Fibrations associ\'ees \`a un pinceau de courbes planes, Annales de la Facult\'e des Sciences de Toulouse, 10 (2001), 745--777.
  
\bibitem{M} J. Milnor, Singular points of complex hypersurfaces.
  Annals of Mathematics Studies 61, Princeton University Press (1968).
    
\bibitem{NZ} A. N\'emethi, A. Zaharia, Milnor fibration at
  infinity.  Indag. Math. 3 (1992), 323--335.
  
\bibitem{Pa} A. Parusi\'nski, On the bifurcation set of complex
  polynomial with isolated singularities at infinity.
  Compositio Math.  97 (1995), 369--384.
  
\bibitem{P} A. Pichon, Real analytic germs $f\bar g$ and
  open-book decompositions of the $3$-sphere, International Journal of Mathematics, 16 (2005),  1--12.  
  
\bibitem{PS} A. Pichon and J. Seade, Fibred multilinks and Singularities $f \bar{g}$. Preprint (2005).
   
\bibitem{Ti} M. Tib\u ar, Regularity at infinity of real and
  complex polynomial functions.  Singularity theory (Liverpool, 1996).
  London Math. Soc. Lecture Note Ser. 263, Cambridge
  Univ. Press (1999), 249--264.

\end{thebibliography}
\end{document}